\documentclass[reqno,12pt]{article}

\usepackage{a4wide}
\usepackage{amsmath,stmaryrd} 
\usepackage{amssymb}
\usepackage{amsthm}
\usepackage[utf8]{inputenc} 
\usepackage{graphicx} 
\usepackage{xcolor}
\usepackage[normalem]{ulem}

\numberwithin{equation}{section}

\newtheorem{theo}{Theorem}

\newtheorem{prop}{Proposition}

\theoremstyle{remark}

\newtheorem*{Remark*}{Remark}
\newtheorem*{Remarks*}{Remarks}

\makeatletter
\newcommand*{\house}[1]{%
 \mathord{%
 \mathpalette\@house{#1}%
 }%
}
\newcommand*{\@house}[2]{%
 \dimen@=\fontdimen8 %
 \ifx#1\scriptscriptstyle\scriptscriptfont
 \else\ifx#1\scriptstyle\scriptfont
 \else\textfont\fi\fi
 3 %
 \sbox0{%
 $#1%
 \vrule width\dimen@\relax
 \overline{%
 \kern2\dimen@
 \begingroup 
 #2%
 \endgroup
 \kern2\dimen@
 }%
 \vrule width\dimen@\relax
 \mathsurround=1.5\dimen@ 
 $%
 }%
 \ht0=\dimexpr\ht0-\dimen@\relax
 \dp0=\dimexpr\dp0+2\dimen@\relax
 \vbox{%
 \kern\dimen@ 
 \copy0 %
 }%
}

\newcommand{\Qbar}{\overline{\mathbb Q}}

\newcommand{\popt}{\lambda}

\begin{document}

\title{A new transcendence measure for the values of the exponential function at algebraic arguments}
\date\today
\author{S. Fischler and T. Rivoal}
\maketitle

\begin{abstract} Let $P\in \mathbb Z[X]\setminus\{0\}$ be of degree $\delta\ge 1$ and usual height $H\ge 1$, and let $\alpha\in \Qbar^*$ be of degree $d\ge 2$. Mahler proved in 1931 the following transcendence measure for $e^\alpha$:  for any $\varepsilon>0$, there exists $c>0$ such that $\vert P(e^\alpha)\vert>c/H^{\mu(d,\delta)+\varepsilon}$ where the exponent 
 $\mu(d,\delta)=(4d^2-2d)\delta+2d-1$. Zheng obtained a better result in 1991 with $\mu(d,\delta)=(4d^2-2d)\delta-1$. In this paper, we provide a new explicit exponent $\mu(d,\delta)$ which improves on Zheng's transcendence measure 
 for all $\delta\ge 2$ and all $d\ge 2$. When $\delta=1$, we recover his bound for all $d\ge 2$, which had in fact already been obtained  by Kappe in 1966.  
Our improvement rests upon the optimization of an accessory parameter in Siegel's classical determinant method applied to Hermite-Padé approximants to powers of the exponential function.  
\end{abstract}
 
\section{Introduction}
In this paper, the field of algebraic numbers $\Qbar$ is embedded into $\mathbb C$. 
Given a number $\xi\in \mathbb C$ known to be transcendental over $\mathbb Q$, a classical way to measure how far $\xi$ is from $\Qbar$ is by giving a non-trivial lower bound of $\vert P(\xi)\vert$ in terms of the height and degree of $P\in\Qbar[X]$. We then say that we have a {\em transcendance measure} of $\xi$. 
There exist many transcendence 
measures for classical numbers such as the values of the exponential function at non-zero algebraic numbers, see for instance \cite{cij, mahler, miw, zheng} and the references given in these papers. It is often difficult to compare these results, because sometimes more attention is paid on the height than on the degree, and vice versa, or also because the measure holds when the height or the degree is assumed large enough with respect to the other. In this paper, we shall provide a  transcendence   measure for $e^{\alpha}$ for all $\alpha\in \Qbar^*$ in which the dependence on the height seems to be the best known so far at this level of generality (Theorem \ref{theo:1} below). 

\medskip

We first introduce some notations. For $\alpha\in \mathbb K$ where $\mathbb K$ is a number field, we define the house of $\alpha$ as $\house{\alpha}=\max_{\sigma}\vert \sigma(\alpha)\vert$ where $\sigma$ runs through all embeddings of $\mathbb K$ into $\mathbb C$; in other words, $\house{\alpha}$ is the maximum of the moduli of $\alpha$ and of all its Galois conjugates over $\mathbb Q$. Given a polynomial $P(X)=\sum_{j=0}^d a_jX^j \in \Qbar[X]$, we set $H(P):=\max_{j=0, \ldots, d}\house{a_j}$ its (usual) height. The (usual) height $H(\alpha)$ of  $\alpha\in \Qbar$ is defined as $H(Q)$ where $Q\in \mathbb Z[X]\setminus\{0\}$ is the minimal polynomial of $\alpha$ over $\mathbb Q$ (normalized so that its coefficients are coprime and the leading coefficient is positive). We also let $\deg(\alpha):=\deg(Q)$ and $d(\alpha)\ge 1$ is the denominator of $\alpha$.

Let $p, d, \delta$ be integers such that $d, \delta,p\ge 1$ and $p\ge \delta d$ and let us define
$$
\psi(d,\delta, p):=\frac{\delta d^2(p-\delta+1)}{p-\delta d+1}+d(p-\delta+1)-1.
$$
This parameter $p$ in the definition of $\psi(d,\delta,p)$ is accessory but appears naturally in the proof of Theorem \ref{theo:1}. It is what enables us to improve on previous bounds. If $d=1$, note that $\psi(1,\delta, p)=p$ is independent of $\delta$, but we still require that $p\ge \delta$.  For $d\ge 2$, we define $p_1:=\delta d-1+\lfloor \delta\sqrt{d^2-d}\rfloor$ and $p_2:=p_1+1$, which are both $\ge \delta d$, and then we set 
$$
\popt :=
\begin{cases} p_1 \quad \textup{if}\quad \psi(d,\delta,p_1)\le \psi(d,\delta, p_2),
\\
p_2 \quad \textup{if}\quad \psi(d,\delta,p_2)<\psi(d,\delta, p_1).
\end{cases}
$$ 
When $d=1$, we set $\popt:=\delta$. We shall prove that $\popt$, which depends on $d$ and $\delta$, minimizes the function $p\mapsto \psi(d,\delta,p)$ amongst all integer values of  $p$ such that $p\ge \delta d$ (this is obvious if $d=1$ by definition).

\medskip

The main result of this article is the following new transcendence measure for $e^\alpha$. 
\begin{theo} \label{theo:1} Let $\mathbb K$ be a number field, let $\alpha\in \Qbar^*$ be such $[\mathbb K(\alpha):\mathbb Q]=d\ge 1$.

For any $\varepsilon>0$ and  any integer $\delta \ge 1$, there exists a constant $c=c(\varepsilon, \alpha, \delta, \mathbb K)>0$ such that for all $H\ge 1$, we have 
\begin{equation}\label{eq:transcmeasure}
  \left\vert P(e^{\alpha})\right\vert  > \frac{c}{H^{\psi(d,\delta,\popt)+\varepsilon}}
\end{equation}
for every polynomial $P\in \mathcal{O}_{\mathbb K}[X]\setminus\{0\}$ of degree $\le \delta$ and height $H(P)\le H$.
\end{theo}
We shall also prove the following facts concerning the exponent $\psi(d,\delta, \lambda)$ in Eq.~\eqref{eq:transcmeasure}. For $d=1$ and all $\delta\ge 1$, $\psi(1,\delta,\popt)=\delta$. For all $d\ge2$ and all $\delta \ge1$, we have 
\begin{multline}\label{eq:boundsf}
\big(2d^2+2d\sqrt{d^2-d}-d\big)\delta-1
\\
\le 
\psi(d,\delta, \popt) \le \big(2d^2+2d\sqrt{d^2-d}-d\big)\delta-1+\frac{d}{\delta \sqrt{d^2-d}-1}.
\end{multline}
Eq.~\eqref{eq:boundsf} will then be used to show that, for all $d\ge 1$, $\psi(d,\delta, \popt)$ is an increasing function of $\delta\ge 1$. 
 The constant $c$ 
 can be made completely explicit, and $\varepsilon$ can be replaced by an explicit decreasing function of $H$; the result being unavoidably complicated, we do not state it in this introduction, and we postpone it to Proposition \ref{prop:1} in \S\ref{sec:effectivisation}. We simply mention here that the dependence in $H$ is of the form $H^{-\psi(d,\delta,\popt)-\varpi/\ln\ln(H+2)}$, where $\varpi>0$ is independent of $H$; such a dependence already occurs in the celebrated transcendence measure of $e$ obtained by Mahler in \cite[p.~135, Satz~3]{mahler}.
 
\medskip

We assume in this paragraph that $\mathbb K=\mathbb Q$ and unless otherwise specified that $\delta\ge 1$ and $d\ge 2$. Another lower bound of the kind in \eqref{eq:transcmeasure}, {\em i.e.}, where the attention is paid on $H$, follows as a special case of the general algebraic independence measure of Lang-Galochkin \cite[p.~238]{feldnest} for $E$-functions: in our situation, it provides the exponent $4\delta d^2$ instead of $\psi(d,\delta,\popt)$ in \eqref{eq:transcmeasure}. In \cite[pp.~132-133]{mahler}, Mahler had obtained in 1931 a smaller exponent, {\em i.e.}, $4\delta d^2-2\delta d+2d-1$ rather than $4\delta d^2$, then improved by Zheng~\cite{zheng} in 1991 with the smaller exponent $(4d^2-2d)\delta -1=:\mu(d,\delta)$. Notice that when $\delta=1$, Kappe \cite{kappe} had already obtained the value $\mu(d,1)$ in 1966~(\footnote{We warn the reader that when $\mathbb K=\mathbb Q$ we measure here $\vert be^\alpha-a\vert$ for $(a,b)\in \mathbb Z^2$, not $\vert e-\frac ab\vert$ as Kappe does, hence the $-1$ in our exponent. Kappe's upper bound also holds for $d=1$ (and it is then an equality to 1) and it is still the best known so far for $d\ge 2$ in general. She obtained it by the different but related method of interpolation series, see also \cite[Chapter II]{schneider}. For certain quadratic numbers like $\sqrt{2}$, Kappe's upper bound 11 has been improved to 3 ({\em i.e.}, $\vert be^{\sqrt{2}}-a\vert>c/b^{3}$) by a different method based on graded Padé approximants to the $E$-function $e^{\sqrt{2}x}+e^{-\sqrt{2}x}$; see \cite{firimu2}.}).
Our exponent  $\psi(d,\delta,\popt)$  is always smaller than or equal to all of these exponents because we shall also prove in \S\ref{ssec:bound4deltadcarre} that $\psi(d,\delta,\popt)\le  (4d^2-2d)\delta-1$ with strict inequality when $\delta\ge 2$ and equality when $\delta=1$. In fact, for all $\delta\ge 2$ and $d\ge 2$, it can be shown that $\psi(d,\delta,\lambda)\le (4d^2-2d-\frac{1}{4})\delta$; this implies that $\mu(d,\delta)-\psi(d,\delta,\lambda)\ge \frac{\delta}{4}-1$ which quantifies more precisely the gain we obtain when $\delta\ge 5$.
If $d=1$, then for all $\delta \ge 1$, the value $\psi(1,\delta,\popt)=\delta$ is the optimal Dirichlet exponent when $\mathbb K=\mathbb Q$, and we recover here a well-known result; see \cite[Chapter 10]{baker}. 

\medskip

Since for all fixed $d\ge 1$, $\psi(d,\delta, \popt)$ is an increasing function of $\delta\ge 1$, it follows that when $\deg(P)\ge 1$ is given, the smallest exponent of $H$ obtainable in \eqref{eq:transcmeasure} is for $\delta=\deg(P)$. In particular, when $d\ge 2$:
\begin{enumerate}
\item[$(i)$] When $\deg(P)=1$, the smallest exponent obtainable in  \eqref{eq:transcmeasure} is $4d^2-2d-1$, obtained for $\delta=1$ with $\popt=2d-2$ (because $\lfloor \sqrt{d^2-d}\rfloor=d-1$ for all $d\ge 2$).
\item[$(ii)$] 
When $\deg(P)=2$, the smallest exponent obtainable in  \eqref{eq:transcmeasure} is 
$$
\frac{16d^3-16d^2+d+1}{2d-1}=8d^2-4d-\frac32+\mathcal{O}\Big(\frac1d\Big),
$$
obtained for $\delta=2$ with $\popt=4d-2$ (because $\lfloor 2\sqrt{d^2-d}\rfloor =2d-2$ for all $d\ge 2$).

\item[$(iii)$] When $\deg(P)=3$, the smallest exponent obtainable in \eqref{eq:transcmeasure} is $$\frac{36d^3-42d^2+7d+2}{3d-1}
=12d^2-6d-\frac53 +\mathcal{O}\Big(\frac1d\Big)
,$$  
obtained for $\delta=3$ with $\popt=6d-3$ (because $\lfloor 3\sqrt{d^2-d}\rfloor =3d-2$ for all $d\ge 2$).
\end{enumerate}
(When $d$ is large enough with respect to $\delta$, we have more generally $\lfloor \delta\sqrt{d^2-d}\rfloor=\delta d -\lceil(\delta+1)/2 \rceil$.) When $\mathbb K=\mathbb Q$, $(i)$ coincides with Kappe's result, but $(ii)$ and $(iii)$ are  better than any previous known bounds.

Finally, it seems that the following properties hold, but we did not try to prove them: for all $d\ge 2$, if $\delta\ge 3$ is odd, $\psi(d,\delta,p_1)<\psi(d,\delta, p_2)$, while if $\delta\ge 2$ is even, $\psi(d,\delta,p_2)<\psi(d,\delta, p_1)$.

\medskip

As already said, Eq.~\eqref{eq:transcmeasure} improves on Zheng's transcendence measure~\cite{zheng} when $\delta\ge 2$ and $d\ge 2$. Like Zheng and Mahler earlier in \cite{mahler}, our proof uses explicit Hermite-Pad\'e type approximants to powers of the exponential function (which have been known explicitly for a long time) and Siegel's classical ``determinant method'' to produce linearly independent linear forms in powers of $e^{\alpha}$. Zheng also considered the parameter $p$ (called $m$ in his paper) which he then took equal to $2\delta d-1$ in his computations. Our contribution lies in the optimization of $p$ with respect to $\delta$ and $d$, and this leads to our improvement. Another interesting application of Siegel's method is in \cite{sorokin}, in which Sorokin computed a very good transcendance measure of $\pi^2$ by this method.

Using a different method ({\em i.e.}, auxiliary functions constructed with Siegel's lemma), Cijsouw \cite[Theorem 1]{cij} obtained a lower bound $\vert P(e^\alpha)\vert >e^{-c(\alpha)\delta^3}H^{-c(\alpha)\delta^2}$ where $P\in \mathbb Z[X]\setminus \{0\}$; his constant $c(\alpha)>0$ is not explicit but from the proof it seems to be larger than $4d^2$ and also to depend on $H(\alpha)$. Our lower bound is thus again better with respect to the dependence on $H$ but, as we shall infer from the examples deduced from Proposition \ref{prop:1} in \S\ref{sec:effectivisation}, his lower bound is better with respect to the dependence on $\delta$. Note that when $H(P)\le e^\delta$, Cijsouw's lower bound has been improved in \cite[p.~455, Corollary~3.9]{miw}.

Besides Zheng's bound, recent works on transcendence measures of $e^\alpha$ essentially all focused on the cases where $\alpha\in \mathbb K$ and $\mathbb K=\mathbb Q$ or a quadratic imaginary field, with results similar to those obtained by Mahler for $e$ in \cite{mahler} or by Baker in \cite{baker2}. Such results improve on Theorem \ref{theo:1} in these particular cases. We refer for instance to \cite{ernvall, ernvall2, hata} and the references therein. 

\medskip

\noindent {\bf Acknowledgement.} We warmly thank Elisabeth Kneller, Yishuai Niu and Arnaud Plessis for providing Zheng's paper \cite{zheng} to us.

\section{Proof of Theorem \ref{theo:1} and related results} \label{sec:proofsthm12}

We first recall basic properties of  Hermite-Pad\'e approximants to powers of the exponential. We then prove Theorem~\ref{theo:1}, and finally  we prove various properties stated in the introduction concerning the exponent $\psi(d,\delta,\lambda)$.

\subsection{Reminder of Hermite-Pad\'e approximants to power of the exponential} \label{ssec:padeapprox}

We state here without proofs properties that are immediate applications of the general construction made in \cite[pp. 63--69]{shid}, applied here with $m:=p+1$ and $\rho_k:=(k-1)\alpha$, $\alpha\in \mathbb C^*$. These properties are necessary for the proof of Theorem \ref{theo:1}.

For all integers $n\ge 1, p\ge 0$, there exist (explicit) polynomials $\mathcal{P}_{k,\ell}\in \mathbb C[x]$ such that for all $\ell \in \{0, \ldots, p\}$
$$
\mathcal{R}_{\ell}(x):=\sum_{k=0}^p \mathcal{P}_{k,\ell}(x) e^{k\alpha x}
$$
vanishes at order $(p+1)n+\ell$ at $x=0$, $\deg(\mathcal{P}_{k,\ell}) = n$ if $k\le \ell$, $\deg(\mathcal{P}_{k,\ell})= n-1$ if $k>\ell$ and 
$$
\det\big(( \mathcal{P}_{k,\ell}(x))_{0\le k,\ell \le p}\big) =cx^{(p+1)n}, \quad c\neq 0.
$$
For simplicity, we do not write explicitly that $\mathcal{R}_{\ell}$ and $\mathcal{P}_{k,\ell}$ also depend on $\alpha$, $n$ and $p$.

\medskip

Now, we assume more specifically that $\alpha\in \mathbb \Qbar^*$, and let $\mathbb K$ be a number field over $\mathbb Q$. We define the positive integer $q:=\textup{lcm}(1,2,\ldots, p)\times d(1/\alpha)$, so that $q^{(p+1)n+p} n! \,\mathcal{P}_{k,\ell}(1)\in \mathcal{O}_{\mathbb Q(\alpha)}\subset  \mathcal{O}_{\mathbb K(\alpha)}$. Then, using the explicit formulas and bounds in \cite[pp.~66-67]{shid}, for all $n\ge 1$, all $p\ge 0$, all $\ell=0, \ldots, p$ and all embeddings $\sigma$ of $\mathbb K(\alpha)$ into $\mathbb C$, we have 
\begin{align} 
 q^{(p+1)n+p} n! \,\sigma(\mathcal{P}_{k,\ell}(1)) &\in \mathcal{O}_{\mathbb K(\alpha)},\label{eq:majoration1}
\\
\vert q^{(p+1)n+p} n! \,\sigma(\mathcal{P}_{k,\ell}(1))\vert &\le \big(2q(1+\vert \sigma(\alpha)\vert^{-1})\big)^{(p+1)n+p} n!,\label{eq:majoration2}
\\
\vert q^{(p+1)n+p} n! \,\mathcal{R}_{\ell}(1)\vert &\le e^{\vert \alpha\vert p(p+1)/2}\cdot \frac{n^{p+1} q^{(p+1)n+p}}{n!^{p}}, \label{eq:majoration3}
\\
\det\big((P_{k,\ell}(1))_{0\le k,\ell \le p}\big) &\neq 0. \label{eq:majoration4}
\end{align}
In \eqref{eq:majoration2}, we can replace $(1+\vert \sigma(\alpha)\vert^{-1})$ by $t_0:=\max_{\sigma}(1+\vert \sigma(\alpha)\vert^{-1})=1+\house{1/\alpha}$ to obtain an upper bound uniform in $\sigma$. Note that $d(1/\sigma(\alpha))$ is independent of $\sigma$, which explain the presence of $q$ in \eqref{eq:majoration1} and \eqref{eq:majoration2}. Note that the upper bound \eqref{eq:majoration2} can be slightly improved to 
\begin{equation} \label{eq:alternative}
\vert q^{(p+1)n+p} n! \,\sigma(\mathcal{P}_{k,\ell}(1))\vert \le (2q\max(1,\vert \sigma(\alpha)\vert^{-1}))^{(p+1)n+p}(n+1)! \le (2qt)^{(p+1)n+p}(n+1)!
\end{equation}
by taking directly $x=1$ in  \cite[p.~67, Eq.~(120)]{shid}, and where the above quantity $t_0$ is now replaced with $t:=\max(1,\house{1/\alpha})$.

Using \cite[p.~80, Lemma~3.11]{miwlivre},  
both quantities $d(1/\alpha)$ and  $\house{1/\alpha}$ can be bounded by $H(\alpha)\sqrt{1+\deg(\alpha)}$, because  $H(\alpha)=H(1/\alpha)$ and $\deg(\alpha)=\deg(1/\alpha)$. 

\subsection{Proof of Theorem \ref{theo:1}} \label{ssec:proofthm11}

Let $\delta \ge 1$ and $p\ge \delta$. (Later on, we shall even impose that $p\ge d\delta$ to obtain our results.) We set $P(X):=\sum_{k=0}^\delta a_k X^k$ with $(a_0,\ldots, a_\delta)\in \mathcal{O}_{\mathbb K}^{\delta+1}\setminus \{0\}$ such that $H(P):=\max\house{a_k}\le H$. Let $\alpha\in \Qbar^*$ and define $d:=[\mathbb K(\alpha):\mathbb Q]$; we have  $L_0:=P(e^\alpha)\neq 0$ because $e^\alpha\notin \Qbar$. 

We also  set $A_{k,\ell}:= q^{(p+1)n+p} n! \,\mathcal{P}_{k,\ell}(1)\in  \mathcal{O}_{\mathbb K(\alpha)} $ and  $R_{\ell}:=q^{(p+1)n+p} n! \,\mathcal{R}_{\ell}(1)$. We have  $\det((A_{k,\ell})_{0\le k,\ell \le p})\neq 0$ because $\det((\mathcal{P}_{k,\ell}(1))_{0\le k,\ell \le p})\neq 0$.

The $p-\delta+1$ vectors  ${}^t(a_0, \ldots, a_\delta, 0,\ldots, 0)$, ${}^t(0, a_0, \ldots, a_\delta, 0,\ldots, 0)$,..., ${}^t(0,\ldots, 0,a_0, \ldots, a_{\delta})$ of $\mathbb C^{p+1}$  are $\mathbb C$-linearly independent. We complete them with the $\delta$ $\mathbb C$-linearly independent vectors  ${}^t(A_{0, \ell}, A_{1, \ell}, \ldots, A_{p,\ell})$ ($\ell=0, \ldots, \delta-1$)  to form a basis of $\mathbb C^{p+1}$.~(\footnote{Strictly speaking, to form a basis of $\mathbb C^{p+1}$, we complete the $p-\delta+1$ first vectors with $\delta$ amongst the $p+1$ vectors ${}^t(A_{0, \ell}, A_{1, \ell}, \ldots, A_{p,\ell})$ ($\ell=0, \ldots, p$). For ease of writing, we take the $\delta$ first ones. In fact, they need not necessarily make up the required basis of $\mathbb C^{p+1}$, but the analysis is completely similar with the $\delta$ good vectors. Moreover, this changes nothing to the effective estimates because Eqs. \eqref{eq:majoration1}-\eqref{eq:majoration4} are uniform in~$\ell$.})

Consequently, the algebraic integer of $\mathbb K(\alpha)$ 
$$
D:=\left\vert 
\begin{matrix}
a_0&a_1&\cdots &a_{\delta}&0&\cdots&\cdots &\cdots&\cdots &0
\\
0 &a_0&\cdots &a_{\delta-1} &a_{\delta}&0&\cdots &\cdots&\cdots &0
\\
\vdots &\vdots&\vdots &\vdots&\vdots&\vdots&\vdots&\vdots&\vdots&\vdots
\\
0 &0&\cdots &0&\cdots&\cdots &0&a_0&\cdots & a_{\delta}
\\
A_{0,0}&A_{1,0}&\cdots &\cdots&\cdots&\cdots &\cdots &\cdots &\cdots& A_{p,0}
\\
A_{0,1}&A_{1,1}&\cdots &\cdots &\cdots \cdots &\cdots &\cdots&\cdots&\cdots& A_{p,1}
\\
\vdots &\vdots&\vdots &\vdots&\vdots&\vdots&\vdots&\vdots&\vdots&\vdots
\\
A_{0,\delta-1}&A_{1,\delta-1}&\cdots &\cdots&\cdots \cdots &\cdots&\cdots &\cdots &\cdots& A_{p,\delta-1}
\end{matrix}
\right\vert \neq 0.
$$
For every embedding $\sigma$ of $\mathbb K(\alpha)$ into $\mathbb C$, we also have $\sigma(D)\neq 0$ so that
$$
\prod_{\sigma} \sigma(D) \in \mathbb Z \setminus \{0\}.
$$
In this product, which is over all such embeddings, we shall distinguish $D$ from the other $\sigma(D)$ with $\sigma\neq id$.

Let $L_j:=\sum_{k=0}^\delta a_k e^{(k+j)\alpha}$. On the one hand, by linear combinations of columns, we find  
$$
D=\left\vert 
\begin{matrix}
L_0&a_1&\cdots &a_{\delta}&0&\cdots&\cdots&\cdots &\cdots &0
\\
L_1 &a_0&\cdots &a_{\delta-1}&a_{\delta}&0&\cdots&\cdots&\cdots &0
\\
\vdots &\vdots&\vdots &\vdots&\vdots&\vdots&\vdots&\vdots&\vdots&\vdots
\\
L_{p-\delta} &0&\cdots &0&\cdots&\cdots&0&a_0&\cdots & a_{\delta}
\\
R_{0}&A_{1,0}&\cdots &\cdots &\cdots&\cdots&\cdots &\cdots&\cdots& A_{p,0}
\\
R_{1}&A_{1,1}&\cdots &\cdots &\cdots&\cdots&\cdots &\cdots&\cdots& A_{p,1}
\\
\vdots &\vdots&\vdots &\vdots&\vdots&\vdots&\vdots&\vdots&\vdots&\vdots
\\
R_{\delta-1}&A_{1,\delta-1}&\cdots &\cdots &\cdots& \cdots&\cdots &\cdots &\cdots& A_{p,\delta-1}
\end{matrix}
\right\vert.
$$
Expanding this determinant about its first column and since $L_j=e^{j\alpha} L_0$, the bounds in Eqs.~\eqref{eq:majoration1}-\eqref{eq:majoration3} imply that 
\begin{equation} \label{eq:majorationD}
\vert D\vert \le c_0^{n+1} H^{p-\delta}n!^\delta \vert L_0\vert + \frac{c_0^{n+1}H^{p-\delta+1}}{n!^{p-\delta+1}}
\end{equation}
for a constant $c_0\ge 1$ independent of $n$ and $H$.

On the other hand, since 
$$
\sigma(D)=\left\vert 
\begin{matrix}
\sigma(a_0)&\sigma(a_1)&\cdots &\sigma(a_{\delta})&0&\cdots&\cdots &\cdots&\cdots &0
\\
0 &\sigma(a_0)&\cdots &\sigma(a_{\delta-1}) &\sigma(a_{\delta})&0&\cdots &\cdots &\cdots &0
\\
\vdots &\vdots&\vdots &\vdots&\vdots&\vdots&\vdots&\vdots&\vdots&\vdots
\\
0 &0&\cdots &0&\cdots&\cdots&0&\sigma(a_0)&\cdots &\sigma(a_{\delta})
\\
\sigma(A_{0,0})&\sigma(A_{1,0})&\cdots &\cdots&\cdots &\cdots &\cdots&\cdots &\cdots& \sigma(A_{p,0})
\\
\sigma(A_{0,1})&\sigma(A_{1,1})&\cdots &\cdots &\cdots &\cdots &\cdots&\cdots&\cdots& \sigma(A_{p,1})
\\
\vdots &\vdots&\vdots &\vdots&\vdots&\vdots&\vdots&\vdots&\vdots&\vdots
\\
\sigma(A_{0,\delta-1})&\sigma(A_{1,\delta-1})&\cdots &\cdots&\cdots&\cdots &\cdots &\cdots &\cdots& \sigma(A_{p,\delta-1})
\end{matrix}
\right\vert, 
$$
we have  
\begin{equation} \label{eq:majorationsigmaD}
\vert \sigma(D)\vert \le c_0^{n+1}H^{p-\delta+1}n!^\delta
\end{equation}
where the constant $c_0\ge 1$ can be taken the same as before (up to  increasing it if necessary). 

Therefore, from $\vert D\vert \prod_{\sigma\neq id} \vert \sigma(D)\vert \ge 1$ and with $d:=[\mathbb K(\alpha):\mathbb Q]$, we deduce that 
$$
\frac{1}{\big(c_0^{n+1}H^{p-\delta+1}n!^\delta\big)^{d-1}} \le \vert D\vert \le c_0^{n+1} H^{p-\delta}n!^\delta \vert L_0\vert + \frac{c_0^{n+1}H^{p-\delta+1}}{n!^{p-\delta+1}}.
$$
Hence, 
\begin{equation}\label{eq:minL}
\vert L_0\vert \ge \frac{1}{c_0^{d(n+1)}H^{(p-\delta+1)(d-1)+p-\delta}n!^{\delta d}}-\frac{H}{n!^{p+1}}=:M
\end{equation}
The right-hand side satisfies
\begin{equation}\label{eq:minM}
M\ge \frac{1}{2c_0^{d(n+1)}H^{(p-\delta+1)(d-1)+p-\delta}n!^{\delta d}}
\end{equation}
provided we can choose $n$ (minimal to get a lower bound as large as possible) such that 
\begin{equation}\label{eq:minN}
    2H^{(p-\delta+1)d}\le n!^{p-d\delta+1} c_0^{-d(n+1)}.
\end{equation}
Since $H\ge 1$ is arbitrary, a necessary condition for the existence of such an $n$ in all circonstances is that $p-\delta d+1>0$, {\em i.e.}, that $p\ge \delta d$ because they are integers. We thus now assume that $p\ge \delta d$ and choose $n$ minimal such that 
\eqref{eq:minN} is satisfied. 
Combining \eqref{eq:minL} and \eqref{eq:minM} with this value of $n$, by standard computations (see \cite[p.~359]{shid} or \S\ref{sec:effectivisation} below), we finally obtain that for all $\varepsilon>0$, there exists a constant $c=c(\varepsilon, \alpha, \delta, \mathbb K)>0$ independent of $H$ such that
$$
\vert L_0 \vert \ge \frac{c}{H^{\psi(d,\delta,p)+\varepsilon}},
$$
where 
$$
\psi(d,\delta, p):=\frac{\delta d^2(p-\delta+1)}{p-\delta d+1}+d(p-\delta+1)-1.
$$
It remains to find the minimal possible value of $\psi(d,\delta, p)$ under the assumption that $p\ge \delta d$. We recall that when $d=1$, $\psi(1,\delta, p)=\delta$ for all $p\ge \delta\ge 1$ so that the minimal value of $\psi(d,\delta, p)$ is achieved at $p=\delta$. We now assume that $d\ge 2$ and $\delta\ge 1$ are fixed. Then the minimum of $x\mapsto \psi(d,\delta,x)$ is attained at 
$$
x_0:=\delta d-1+\delta\sqrt{d^2-d},
$$
and the integers $p_1:=\lfloor x_0\rfloor=\delta d-1+\lfloor \delta\sqrt{d^2-d}\rfloor$ and $p_2:=\lfloor x_0\rfloor+1=\delta d+\lfloor \delta\sqrt{d^2-d}\rfloor$ are both admissible to minimize $\psi(d,\delta,x)$ with $x$ an integer (because both are $\ge \delta d$). Therefore defining $\popt$ as either $p_1$ if $\psi(d,\delta,p_1)\le \psi(d,\delta, p_2)$ or $p_2$ if $\psi(d,\delta,p_2)<\psi(d,\delta, p_1)$, we obtain that $\psi(d,\delta, p)\ge \psi(d,\delta,\popt)$ for all $p\ge \delta d$.  This completes the proof of Theorem~\ref{theo:1}.

\subsection{The upper bound in \eqref{eq:boundsf}} \label{ssec:proofthm12}

Direct computations show that this is true for all $d\ge 2$ if $\delta=1,2$, because 
$\psi(d,1,\popt)=4d^2-2d-1$ and 
$$
\psi(d,2,\popt)=\frac{16d^3-16d^2+d+1}{2d-1}
$$
(because $\lfloor \sqrt{d^2-d}\rfloor=d-1$ 
and $\lfloor 2\sqrt{d^2-d}\rfloor=2d-2$), 
which are both $\le \big(2d^2+2d\sqrt{d^2-d}-d\big)\delta-1+\frac{d}{\delta \sqrt{d^2-d}-1}$ for $\delta=1,2$ and $d\ge 2$. 
We now assume that $\delta \ge 3$ and $d\ge 2$.

\medskip

We start with $p_1$. We write 
$$
p_1=\delta d +\delta \sqrt{d^2-d}+x, \quad x\in(-1,0]
$$
and $D:=\sqrt{d^2-d}$ to simplify.
Then 
\begin{align}
\psi(d,\delta, p_1)&=\frac{\delta d^2(\delta d-\delta+\delta D+x)}{\delta D+x}+d(\delta d+\delta D+x-\delta)-1 \notag
\\
&=2\delta d^2+\delta dD-\delta d -1+ \frac{\delta^2d^2(d-1)}{\delta D+x}+dx. \label{eq:minorationfp1}
\end{align}
The function $x\mapsto \frac{\delta^2d^2(d-1)}{\delta D+x}+dx$ is decreasing on the interval $[-1,0]$ (its derivative is $\frac{dx(2\delta D+x)}{(\delta D+x)^2}$) so that 
$$
\psi(d,\delta, p_1) \le 2\delta d^2+\delta d\sqrt{d^2-d}-\delta d -1+ \frac{\delta^2d^2(d-1)}{\delta \sqrt{d^2-d}-1}-d.
$$
Now 
$$
\frac{\delta^2d^2(d-1)}{\delta \sqrt{d^2-d}-1} = \delta d\sqrt{d^2-d}+d+\frac{d}{\delta\sqrt{d^2-d}-1}
$$
so that 
\begin{equation} \label{eq:majpsi1}
\psi(d,\delta, p_1)\le \big(2 d^2+2 d\sqrt{d^2-d} - d\big)\delta-1 +\frac{d}{\delta\sqrt{d^2-d}-1} 
\end{equation}

\medskip

For $p_2=p_1+1$, we write
$$
p_2=\delta d +\delta \sqrt{d^2-d}+x, \quad x\in(0,1].
$$
Then again 
\begin{align}
\psi(d,\delta, p_2)&=\frac{\delta d^2(\delta d-\delta+\delta D+x)}{\delta D+x}+d(\delta d+\delta D+x-\delta)-1 \notag
\\
&=2\delta d^2+\delta dD-\delta d -1+ \frac{\delta^2d^2(d-1)}{\delta D+x}+dx. \label{eq:minorationfp2}
\end{align}
Since the function $x\mapsto \frac{\delta^2d^2(d-1)}{\delta D+x}+dx$ is now increasing on the interval $[0,1]$, we have 
$$
\psi(d,\delta, p_2) \le 2\delta d^2+\delta d\sqrt{d^2-d}-\delta d-1 + \frac{\delta^2d^2(d-1)}{\delta \sqrt{d^2-d}+1}+d.
$$
Now 
$$
\frac{\delta^2d^2(d-1)}{\delta \sqrt{d^2-d}+1} = \delta d\sqrt{d^2-d}-d+\frac{d}{\delta\sqrt{d^2-d}+1}
$$
so that 
\begin{equation}\label{eq:majpsi2}
\psi(d,\delta, p_2)\le \big(2d^2+2d\sqrt{d^2-d} -d)\delta-1 +\frac{d}{\delta\sqrt{d^2-d}+1}.
\end{equation}
We remark that the right-hand side of \eqref{eq:majpsi1} is larger than the right-hand side of \eqref{eq:majpsi2}, so that 
$$
\psi(d,\delta, \popt)\le \big(2d^2+2d\sqrt{d^2-d} -d\big)\delta-1 +\frac{d}{\delta\sqrt{d^2-d}-1} 
$$
as claimed.

\subsection{The lower bound in \eqref{eq:boundsf}} \label{ssec:proofthm13}

If $d\ge 2$ and $\delta=1,2$, the values $\psi(d,\delta, \popt)$ recalled in \S\ref{ssec:proofthm12} (for $\popt=2d-2$ and $\popt=4d-2$ respectively) are both  $\ge (2d^2+2d\sqrt{d^2-d}-d)\delta-1$. For $d\ge2$ and $\delta\ge3$, we observe that in \eqref{eq:minorationfp1} and \eqref{eq:minorationfp2}, we obtain lower bounds for 
$\psi(d,\delta, p_1)$ and $\psi(d,\delta, p_2)$  respectively by taking $x=0$. It turns out that these lower bounds are both equal to $(2d^2+2d\sqrt{d^2-d}-d)\delta-1$.  

\subsection{$\delta\mapsto \psi(d,\delta, \popt)$ is increasing} \label{ssec:psiincreasing}

We first recall that $\popt$ depends on $d$ and $\delta$, and below we write $\popt_\delta$ for $\popt$. Since $\psi(1,\delta, \popt_\delta)=\delta$, the property holds for $d=1$. Let us now assume that $d\ge 2$. We shall prove that for all $\delta\ge 1$, we have $\psi(d,\delta+1, \popt_{\delta+1})>\psi(d,\delta, \popt_\delta)$. For this, it is enough to prove that the left-hand side of \eqref{eq:boundsf} with $\delta$ replaced by $\delta+1$ is larger than the right-hand side of \eqref{eq:boundsf} with $\delta$. To prove this, we simply observe that
\begin{align*}
\big(2d^2+2d\sqrt{d^2-d}-d&\big)(\delta+1)-\big(2d^2+2d\sqrt{d^2-d}-d\big)\delta-\frac{d}{\delta \sqrt{d^2-d}-1}\\
&=2d^2+2d\sqrt{d^2-d}-d - \frac{d}{\delta \sqrt{d^2-d}-1}
\\
&\ge 2d^2+2d\sqrt{d^2-d}-d - \frac{d}{\sqrt{d^2-d}-1} >0
\end{align*}
for all $d\ge2$ and all $\delta\ge 1$.

\subsection{The inequality $\psi(d,\delta,\popt)\le 4\delta d^2-2\delta d-1$} \label{ssec:bound4deltadcarre} 
 
If $d=1$, $\psi(1,\delta,\popt)=\delta \le 2\delta-1$ for all $\delta \ge 1$.  If $d\ge 2$, we simply observe that $\psi(d,\delta, \popt)\le \psi(d,\delta,2\delta d-1)=4\delta d^2-2\delta d-1$. The inequality is an equality when $\delta=1$ and is strict if $\delta\ge 2$ because then $\lambda\le \delta d+\lfloor \delta \sqrt{d^2-d}\rfloor< 2\delta d-1$. Notice that $p:=2\delta d-1$ is exactly the choice made by Zheng.

\section{A completely explicit version of Theorem \ref{theo:1}} \label{sec:effectivisation}

\subsection{The main statement}
The goal of this section is to make completely explicit the constant $c$ in Theorem~\ref{theo:1}. We shall prove the following proposition, which, in fact, subsumesmess Theorem \ref{theo:1}.
\begin{prop} \label{prop:1}
Let $\mathbb K$ be a number field, let $\alpha\in \Qbar^*$ be such $[\mathbb K(\alpha):\mathbb Q]=:d\ge 1$.

For any $H\ge 1$, any integer $\delta\ge 1$, and integer $p$ such that $p\ge d\delta$, and any $P\in \mathcal{O}_{\mathbb K}[X]\setminus\{0\}$ of degree $\le \delta$ and height $H(P)\le H$, we have 
\begin{equation}
\left\vert P(e^{\alpha})\right\vert  \ge \frac{1}{(2a^d)^{1+\frac{\delta d}{p-d\delta+1}}\big(b^{d+\frac{\delta d^2}{p-d\delta+1}}\big)^{1+v^2} u^{(d+\frac{\delta d^2}{p-d\delta+1})\frac{4\ln(b)}{\ln\ln(u+2)}}H^{\psi(d,\delta,p)}}, \label{eq:mesL0effective*}
\end{equation}
where  
\begin{align*}
q&:=\textup{lcm}(1,2,\ldots, p)\times d(1/\alpha), \qquad 
t:=\max(1,\house{1/\alpha})
\\
a&:=(p+1)!e^{\vert \alpha\vert p(p+1)/2}(2qt)^{p\delta}, \qquad 
b:=(2qt)^{(p+1)\delta},
\\
u&:=\left(2a^dH^{(p-\delta+1)d}\right)^{1/(p-d\delta+1)}, \qquad 
v:=b^{d/(p-d\delta+1)}e,
\\
\psi&(d,\delta,p):=\frac{\delta d^2(p-\delta+1)}{p-\delta d+1}+d(p-\delta+1)-1.
\end{align*}
\end{prop}
After simplifications, the dependence on $H$ turns out to be of the standard ``Mahlerian'' form $$H^{-\psi(d,\delta, p)-\varpi/\ln\ln(H+2)}$$ for some $\varpi>0$ independent of $H$. In \cite[Theorem~1]{zheng}, when $\mathbb K=\mathbb Q$, Zheng has obtained a lower bound of the form 
$H^{-(4\delta d^2-2\delta d-1)-c\delta^2/\ln\ln(H+2)}$, where $c>0$ is an unspecified constant. We recall that, with respect to the accessory parameter $p$, the minimal value of $p\mapsto\psi(d,\delta, p)$ is obtained at $p_1:=\delta d-1+\lfloor \delta\sqrt{d^2-d}\rfloor$ or $p_2=\delta d+\lfloor \delta\sqrt{d^2-d}\rfloor$.

\begin{proof} The proof follows the same steps as that of Theorem \ref{theo:1}, except that  we pay attention to effectivity and explicit bounds. For this we shall use the explicit bounds 
recalled in Eqs.~\eqref{eq:majoration1}-\eqref{eq:majoration3} in \S\ref{ssec:padeapprox}. We can then replace \eqref{eq:majorationD} by the completely explicit bound, valid for all $n\ge 1, p\ge \delta \ge 1, H\ge 1$:
\begin{equation*}
\vert D\vert \le c_1 c_2^{n} H^{p-\delta}n!^\delta \vert L_0\vert + \frac{c_3c_4^{n}n^{p+1}H^{p-\delta+1}}{n!^{p-\delta+1}},
\end{equation*}
where 
\begin{align*}
c_1&:=(p+1)!e^{\vert \alpha\vert p}(2qt)^{p\delta}, \quad
c_2:=(2qt)^{(p+1)\delta}, 
\\
c_3&:=(p+1)! e^{\vert \alpha\vert p(p+1)/2}q^{p}(2qt)^{p(\delta-1)}, \quad 
c_4:=q^{p+1}(2qt)^{(p+1)(\delta-1)}.
\end{align*}
Similarly, we have the explicit bound
\begin{equation*} 
\vert \sigma(D)\vert \le c_5 c_2^{n}H^{p-\delta+1}n!^\delta,
\end{equation*}
for all embeddings $\sigma$ of $\mathbb K(\alpha)$ into $\mathbb C$, where $c_5=(p+1)!(2qt)^{p\delta}$.

We now make a few simplifications. We have $n^{p+1}< (2t)^{(p+1)n}$ for all $n\ge 1$ so that $n^{p+1}c_4^{n}<((2t)^{p+1}c_4)^{n}=c_2^n$. Moreover, $\max(c_1,c_3, c_5)<c_6$ where
$$
c_6:=(p+1)!e^{\vert \alpha\vert p(p+1)/2}(2qt)^{p\delta}.
$$
Therefore,
\begin{equation} \label{eq:majorationDbis}
\vert D\vert \le c_6 c_2^{n} H^{p-\delta}n!^\delta \left(\vert L_0\vert + \frac{H}{n!^{p+1}}\right)
\end{equation}
and 
\begin{equation*} \label{eq:majorationsigmaDbis}
\vert \sigma(D)\vert \le c_6 c_2^{n}H^{p-\delta+1}n!^\delta
\end{equation*}
for all embeddings $\sigma$ of $\mathbb K(\alpha)$ into $\mathbb C$. 
For simplicity, we now set $b:=c_2$ and $a:=c_6$.

From the lower bound $\vert D\vert \prod_{\sigma\neq id} \vert \sigma(D)\vert \ge 1$, where the product is over all such embeddings, we deduce that 
$$
\frac{1}{\big(ab^{n}H^{p-\delta+1}n!^\delta\big)^{d-1}} \le \vert D\vert \le 
a b^{n} H^{p-\delta}n!^\delta \left(\vert L_0\vert + \frac{H}{n!^{p+1}}\right)
$$
where $d:=[\mathbb K(\alpha):\mathbb Q]$. 
Hence, 
\begin{equation}\label{eq:minLbis}
\vert L_0\vert \ge \frac{1}{a^d b^{dn}H^{(p-\delta+1)(d-1)+p-\delta}n!^{\delta d}}-\frac{H}{n!^{p+1}}=:M
\end{equation}
The right-hand side of \eqref{eq:minLbis} satisfies
\begin{equation}\label{eq:minMbis}
M\ge \frac{1}{2 a^d b^{dn}H^{(p-\delta+1)(d-1)+p-\delta}n!^{\delta d}}
\end{equation}
provided we can choose $n$ (as small as possible to get a lower bound as large as possible) such that 
\begin{equation}\label{eq:minNbis}
    2a^dH^{(p-\delta+1)d}\le n!^{p-d\delta+1} b^{-dn}.
\end{equation}
Since $H\ge 1$ is arbitrary, a necessary condition for the existence of such an $n$ in all circumstances is that $p-d\delta+1>0$, {\em i.e.}, that $p\ge d\delta$ because $d,\delta, p$ are integers. We thus now assume that $p\ge d\delta$ and ideally we choose the minimal value of $n$, say $n_{\min}$, such that \eqref{eq:minNbis} is satisfied. Finding the exact expression on $n_{\min}$ is difficult but we can find an upper bound for $n_{\min}$ as follows. Since $n!\ge n^n e^{-n}$ for all $n\ge 1$, the minimal value of $n$ (denoted by $n_0$ from now on) such that 
\begin{equation*}
    2a^dH^{(p-\delta+1)d}\le n^{(p-d\delta+1)n} e^{-(p-d\delta+1)n}b^{-dn}
\end{equation*}
is such that $n_{\min}\le n_0$. Hence \eqref{eq:minNbis} holds with $n=n_0$, and thus \eqref{eq:minMbis} as well.

Our task is now to obtain an upper bound for $n_0$. By definition of $n_0$, we have
$$
(n_0-1)^{n_0-1}<u v^{n_0-1}
$$
where 
$$
u:=\left(2a^dH^{(p-\delta+1)d}\right)^{1/(p-d\delta+1)}>1, \quad v:=b^{d/(p-d\delta+1)}e>e.
$$
The integer $n_0$ is obviously $\ge v+1$ because $v^v<uv^v$. We then set $n=v+x+1$ for some $x\ge 1$ so that for all $s\in (0,1)$, we have:
\begin{align*}
\frac{(n-1)^{n-1}}{u v^{n-1}} &=\frac{(v+x)^{v+x}}{u v^{v+x}}
\\
&\ge \frac{(v+x)^{(v+x)(1-s)}}{u} \qquad \textup{provided}\; v+x\ge v^{1/s}
\\
&\ge \left(\frac{x^{x}}{u^{1/(1-s)}}\right)^{1-s} \qquad \textup{because}\; x\ge 1/e
\\
&\ge 1 \qquad \textup{provided}\;  x\ge 1+ \frac{2\ln(u^{1/(1-s)})}{\ln\ln(u^{1/(1-s)}+2)}.
\end{align*}
In the last two lines, we use two elementary facts: 1) the function $v\mapsto (v+x)^{v+x}$ increases for $[0,+\infty)$ when $x\ge 1/e$, and 2) for any $h>1$, if $x\ge 1+ \frac{2\ln(h)}{\ln\ln(h+2)}$, then $x^x\ge h$. The three assumptions on $x$ are then fulfilled with 
$$
x_s:=v^{1/s}-v+\frac{2\ln(u)}{(1-s)\ln\ln(u^{1/(1-s)}+2)}+1.
$$
It follows that for all $s\in (0,1)$, we have  
$$
n_0 \le v+x_s = v^{1/s}+ \frac{2\ln(u)}{(1-s)\ln\ln(u^{1/(1-s)}+2)}+1.
$$ 
Taking $s=1/2$, we deduce the upper bound:
\begin{equation}\label{eq:boundn0}
n_0 \le v^2+ \frac{4\ln(u)}{\ln\ln(u^2+2)}+1\le  v^2+ \frac{4\ln(u)}{\ln\ln(u+2)}+1.
\end{equation}

Since \eqref{eq:minNbis} is satisfied with $n=n_0$, we have 
$$
n_0!\ge (2a^dH^{(p-\delta+1)d})^{1/(p-d\delta+1)}b^{d n_0/(p-d\delta+1)}.
$$
Substituting this into \eqref{eq:minMbis}, we deduce from \eqref{eq:minLbis} that
\begin{equation*}
\vert L_0\vert \ge \frac{1}{(2a^d)^{1+\delta d/(p-d\delta+1)}(b^{d+\delta d^2/(p-d\delta+1)})^{n_0}H^{\psi(d,\delta,p)}}. 
\end{equation*}
Using the upper bound for $n_0$ in \eqref{eq:boundn0}, we then obtain
\begin{equation*}
\vert L_0\vert \ge \frac{1}{\left(2a^d\right)^{1+\frac{\delta d}{p-d\delta+1}}\big(b^{d+\frac{\delta d^2}{p-d\delta+1}}\big)^{1+v^2} u^{(d+\frac{\delta d^2}{p-d\delta+1})\frac{4\ln(b)}{\ln\ln(u+2)}}H^{\psi(d,\delta,p)}}, 
\end{equation*}
which is the expected lower bound \eqref{eq:mesL0effective*}. 
\end{proof}

\subsection{Some consequences of Proposition \ref{prop:1}} \label{ssec:examples}
In what follows, $\mathbb K=\mathbb Q$. We present without details three examples of application of Proposition \ref{prop:1}.

\medskip

a) If $\alpha=1$, we have $d=1$, $p=\delta$ and $\psi(d,\delta,p)=\delta$, $t=2$, $q=\textup{lcm}(1,2,\ldots, \delta)\approx e^\delta$, $a \approx e^{\delta^3}$, $b \approx e^{\delta^3}$, $u\approx H^d e^{\delta^3}$, $v\approx e^{d\delta^3}$. We deduce that there exist two absolute constants $c_1>0, c_2>0$
such that for any $H\ge1$, any integer $\delta \ge 1$ and any vector $(a_0, \ldots, a_\delta)\in \mathbb Z^{\delta+1}\setminus\{0\}$ with $\max\vert a_k\vert\le H$, we have 
$$
\left\vert \sum_{k=0}^\delta a_k e^{k} \right\vert \ge \frac{\exp(-
\exp(c_1\delta^3))}{H^{\delta+\frac{c_2\delta^5}{\ln\ln(H+2)}}}.
$$
A completely similar result holds for any $\alpha\in \mathbb Q^*$, except that $c_1$ and $c_2$ now depend on $\alpha$. 
Note that Mahler \cite[p. 135, Satz 3]{mahler} obtained a better exponent when $\alpha=1$: 
$$
\left\vert \sum_{k=0}^\delta a_k e^{k} \right\vert \ge \frac{1}{H^{\delta+\frac{c_3\delta^2\ln(\delta+1)}{\ln\ln(H+2)}}}
$$
where $c_3>0$ is absolute and $H\ge H(\delta)$ (which is not explicit). 

\medskip 

b) We fix $\alpha$ of degree $d\ge 2$ over $\mathbb Q$. In this case, $p\approx 2 d\delta$, $q\approx ce^{2 d\delta}$, $t\approx c$, $a\approx e^{cd^2\delta^3}$, $b\approx e^{cd^2\delta^3}$, $u\approx e^{cd^2\delta^2}H^{2d}$, $v\approx e^{cd^2\delta^2}$ (where $c$ denotes a constant that depends on $\alpha$ but not on $d$, and may differ in each of the previous estimates). It follows that there exist two constants $c_4, c_5>0$, that depend on $\alpha$ but not on $d$, such that for any $H\ge1$, any integer $\delta \ge 1$ and any vector $(a_0, \ldots, a_\delta)\in \mathbb Z^{\delta+1}\setminus\{0\}$ with $\max\vert a_k\vert\le H$, we have
$$
\left \vert \sum_{k=0}^\delta a_k e^{k\alpha}\right\vert \ge \frac{\exp(- \exp (c_4d^2\delta^2))}{H^{\psi(d,\delta,\popt)+\frac{c_5d^4\delta^3}{\ln\ln(H+2)}}}.
$$

\medskip

c) If $\delta=1$, then $p=2d-2$ and $\psi(d,\delta, 2d-2)=4d^2-2d-1$. Estimating the parameters as above and using the fact that $d(1/\alpha)$ and $\house{1/\alpha}$ are both bounded by $\sqrt{d+1}H(\alpha)$, we deduce the existence of two  absolute constants $c_6, c_7>0$
such that for all $\alpha\in \Qbar^*$ of height 
$H(\alpha)$ and degree $d$, for all $(p,q)\in \mathbb Z\times \mathbb N^*$, we have 
$$
\left\vert e^\alpha- \frac{p}{q} \right\vert \ge \frac{\exp\big(- 
\exp(c_6 s(\alpha)d)\big)}{q^{4d^2-2d+\frac{c_7d^3s(\alpha)}{\ln\ln(q+2)}}},
$$
where $s(\alpha):=d+\ln H(\alpha)$ is the classical quantity called the size of $\alpha$ (see \cite{cij}).
This is an explicit version of Kappe's irrationality measure \cite{kappe}.

\section{A more general construction}

We conclude this paper with a more general construction of $\mathbb K$-linear approximations to powers of $e^\alpha$. However, we shall prove that it does not lead to better results than those obtained in the previous sections.

We start as in \S\ref{sec:proofsthm12}, of which we borrow the notations. For any integer $1\le k\le p-\delta+1$,  the $k$ vectors ${}^t(a_0, \ldots, a_\delta, 0,\ldots, 0)$, ${}^t(0, a_0, \ldots, a_\delta, 0,\ldots, 0)$, ..., ${}^t(0,\ldots, 0,a_0, \ldots, a_\delta,0,\ldots, 0)$  of $\mathbb{C}^{p+1}$ are $\mathbb C$-linearly independent. We complete them with the $p-k+1$ $\mathbb C$-linearly independent vectors ${}^t(A_{0, \ell}, A_{1, \ell}, \ldots, A_{p,\ell})$ ($\ell=0, \ldots, p-k$)   of $\mathbb{C}^{p+1}$ to form a basis of $\mathbb C^{p+1}$. The previous construction corresponds to the case $k=p-\delta+1$. 

The integers $p,d,\delta,k$ must satisfy the conditions: 
$$
d\ge1, \quad p\ge \delta\ge 1, \quad 1\le k\le p-\delta+1, \quad p+kd-d(p+1)+1>0.
$$
We no longer require that $p\ge d\delta$.

Then similar ``determinantal'' computations as in the proof of Theorem \ref{theo:1} show that 
$
\vert L_0\vert \ge cH^{-\varphi(d,\delta,p,k)-\varepsilon}
$
where 
$$
\varphi(d,\delta,p,k):=dk-1 + \frac{d^2k(p-k+1)}{p+kd-d(p+1)+1}.
$$
The parameter $\delta$ does not appear explicitly in $\varphi$, and only through the above inequalities. 

If $d=1$, $\varphi(d,\delta,p,k)=p$ is minimal for $p=\delta$ (as expected) and it is enough to take $k=1$ to obtain this result in this case. 

If $d\ge2$, we have
$$
\frac{\partial \varphi}{\partial k}= - \frac{d(d-1)(p+1)^2}{(p+kd-d(p+1)+1)^2}<0.
$$
Therefore, for fixed $d\ge2$ and $p, \delta\ge 1$, $\varphi(d,\delta,p,k)$ is minimal for $k=p-\delta+1$ in which case it coincides with $\psi(d,\delta,p)$.  

In both cases, this justifies to consider only the case $k=p-\delta+1$ in the previous sections.

\noindent St\'ephane Fischler, Universit\'e Paris-Saclay, CNRS, Laboratoire de math\'ematiques d'Orsay, 91405 Orsay, France.

\medskip

\noindent Tanguy Rivoal, Universit\'e Grenoble Alpes, CNRS, Institut Fourier, CS 40700, 38058 Grenoble cedex 9, France.

\bigskip

\noindent Keywords: Exponential function,  Transcendence measure, Hermite-Pad\'e approximants.

\bigskip

\noindent MSC 2020: 11J82 (Primary), 11J91 (Secondary)

\end{document}